\documentclass[journal]{IEEEtran}
\usepackage{amsfonts}
\usepackage{amsmath,amsfonts,amsthm,amssymb}
\usepackage{threeparttable}
\usepackage{bm}
\usepackage{enumitem}
\usepackage{tabularx}
\usepackage{bbm}
\usepackage{float}
\usepackage{setspace}
\usepackage{Tabbing}
\usepackage{fancyhdr}
\usepackage{lastpage}
\usepackage{extramarks}
\usepackage{chngpage}
\usepackage[ruled,linesnumbered]{algorithm2e}
\usepackage{algpseudocode}
\usepackage{amsmath}
\usepackage{soul,color}
\usepackage{epsfig}
\usepackage{graphicx,float,wrapfig,subfigure}
\usepackage{dsfont}
\usepackage{longtable}
\usepackage{wrapfig}
\usepackage{stfloats}
\usepackage{cite}
\usepackage{graphicx}
\usepackage{array}
\usepackage{multirow}{\tiny }
\usepackage{multicol}
\usepackage{colortbl}
\usepackage{tabularx}
\usepackage{mdwmath}
\usepackage{mdwtab}
\usepackage{color}
\usepackage{verbatim}
\usepackage{tikz}
\usetikzlibrary{calc}
\usepackage{flowchart}
\usepackage{algpseudocode}
\usetikzlibrary{shapes.geometric, arrows}
\usepackage{overpic}
\usepackage{epstopdf}
\usepackage{url}
\usepackage{graphicx}
\usepackage{float}
\usepackage{subfigure}
\newcommand{\tabincell}[2]{\begin{tabular}{@{}#1@{}}#2\end{tabular}}
\usepackage[colorlinks,linkcolor=black,anchorcolor=black,citecolor=black,urlcolor=black]{hyperref}
\newtheorem{proposition}{Proposition}

\usepackage{color}
\makeatletter % 
\let\myorg@bibitem\bibitem
\def\bibitem#1#2\par{%
	\@ifundefined{bibitem@#1}{%
		\myorg@bibitem{#1}#2\par
	}{%
		\begingroup
		\color{\csname bibitem@#1\endcsname}%
		\myorg@bibitem{#1}#2\par
		\endgroup
	}%
}

\makeatother % 
\allowdisplaybreaks

\usepackage{footnote}
\makesavenoteenv{tabular}
\makesavenoteenv{table}

\SetKwInOut{Input}{Input}
\SetKwInOut{Output}{Output}

\begin{document}

\title{
%	\vspace{-4mm}
Economic Bidding Strategy of Electric Vehicles in Real-Time Electricity Markets based on Marginal Opportunity Value
	}

\author{
Zhen~Zhu,~\IEEEmembership{Graduate Student Member,~IEEE,}
Hongcai~Zhang,~\IEEEmembership{Senior Member,~IEEE}, and
Yonghua~Song,~\IEEEmembership{Fellow,~IEEE}
%\vspace{-4mm}

\thanks{	
Z. Zhu, H. Zhang and Y. Song are with the State Key Laboratory of Internet of Things for Smart City and Department of Electrical and Computer Engineering, University of Macau, Macao, 999078 China (email: hczhang@um.edu.mo).
}
}

\maketitle

\begin{abstract}  
The participation of electric vehicle (EV) aggregators in real-time electricity markets offers promising revenue opportunities through price-responsive energy arbitrage. A central challenge in economic bidding lies in quantifying the marginal opportunity value of EVs' charging and discharging decisions. This value is implicitly defined and dynamically shaped by uncertainties in electricity prices and availability of EV resources. In this paper, we propose an efficient bidding strategy that enables EV aggregators to generate market-compliant bids based on the underlying marginal value of energy. The approach first formulates the EV aggregator's power scheduling problem as a Markov decision process, linking the opportunity value of energy to the value function. Building on this formulation, we derive the probability distributions of marginal opportunity values across EVs' different energy states under stochastic electricity prices. These are then used to construct closed-form expressions for marginal charging values and discharging costs under both risk-neutral and risk-averse preferences. 
% Based on this formulation, we derive the probability distributions of marginal opportunity values for EVA under stochastic prices and develop closed-form expressions for marginal charging values and discharging costs under both risk-neutral and risk-averse preferences. 
The resulting expressions support a fully analytical bid construction procedure that transforms marginal valuations into stepwise price–quantity bids without redundant computation. Case studies using real-world EV charging data and market prices demonstrate the effectiveness and adaptability of the proposed strategy.
\end{abstract}
%\vspace{-4mm}
\begin{IEEEkeywords}
 Electric vehicle aggregator, dynamic risk measure, economic bidding, marginal opportunity value.
\end{IEEEkeywords}

\section{Introduction}
\IEEEPARstart {T}he accelerating transition toward electrified transportation is reshaping the global power system landscape. With millions of electric vehicles (EVs) expected on the road in the coming decades, their collective charging behavior introduces both unprecedented challenges and valuable opportunities for power grid management \cite{review}. In this context, EV aggregators (EVAs), which coordinate large-scale EV fleets, are emerging as critical entities for integrating distributed EV resources into grid operations.
% EV aggregators (EVAs)—entities that coordinate large fleets of EVs—are emerging as key intermediaries, facilitating the integration of decentralized EV resources into centralized grid operations. 
To fully realize the potential of EVs as flexible resources, policymakers worldwide are promoting EVA participation in electricity markets. In the United States, FERC Order No. 2222 allows EVAs to bid directly into wholesale energy markets \cite{ferc}. In China, policies led by the National Development and Reform Commission encourage EVA involvement in spot market transactions \cite{ndrc}. These regulatory efforts not only create new revenue opportunities for EVAs but also support improved power system efficiency and reliability. Consequently, the development of economically sound bidding strategies is becoming increasingly essential for EVA to effectively manage and monetize EV flexibility in market environments.

Energy markets typically consist of day-ahead and real-time segments \cite{wuchen}, both offering opportunities for EVA to generate profits through energy arbitrage. By leveraging temporal price variations, an EVA can shift EV charging to low-price periods and, where vehicle-to-grid (V2G) capabilities are available, discharge energy back to the grid during price peaks, thereby reducing energy costs while fulfilling mobility needs. Recent studies have extensively explored EVA bidding strategies in the day-ahead market, where hourly bids are submitted one day in advance based on forecasts of prices and vehicle availability. For example, reference \cite{dh1} introduces a transactive-based day-ahead scheduling framework that coordinates EV charging using distribution locational marginal prices and EV response curves. Reference \cite{dh2} proposes a risk-averse bidding model that enables EVA to minimize worst-case regret under uncertainties in prices, demand, and renewable generation. Reference \cite{dh3} presents a stochastic day-ahead bidding scheme for a fast-charging station aggregator that incorporates traffic dynamics and bounded EV user rationality. While the day-ahead market provides a relatively stable price environment, it inherently restricts EVA’s ability to respond to real-time price fluctuations, limiting their full arbitrage potential. In contrast, the real-time market is characterized by higher price volatility and shorter settlement intervals, which better align with the fast response capabilities of EV batteries. This enhances arbitrage opportunities and makes real-time market participation increasingly attractive for EVA.
% discharge energy during peak-price intervals, 
% While day-ahead participation provides a stable bidding environment, it inherently limits responsiveness to real-time price fluctuations, restricting the EVA's ability to fully utilize short-term arbitrage opportunities. In contrast, the real-time market features higher price volatility and shorter settlement intervals, which better align with the fast-responding characteristics of EV batteries. These characteristics offer greater arbitrage potential, making real-time market participation increasingly attractive for EVAs.

Unlike the 24-hour commitment structure of day-ahead markets, real-time electricity markets require participants to submit bids sequentially, often just prior to each market interval. This rolling structure enables greater responsiveness to actual system conditions but also introduces unique challenges for the EVA. A primary difficulty stems from the mismatch between the time-coupled nature of EVA operations and the myopic structure of the real-time market, where each bid pertains to only a single interval. Coordinating decisions across time while satisfying inter-temporal constraints such as EV state-of-charge dynamics is inherently complex. Moreover, uncertainties in both real-time electricity prices and EV availability further complicate the bidding process, necessitating strategies that can effectively balance flexibility, foresight, and robustness. Actually, research on EVA bidding strategies in the real-time market remains relatively sparse. Reference \cite{rt3} presents a hierarchical coordination framework, with upper-level market bidding and lower-level power allocation based on a trade-off between aggregator profit and user satisfaction.
In \cite{rt2}, an endpoint energy and power boundary model is introduced to handle cross-day charging discontinuities, enhancing bidding continuity. Reference \cite{rt1} develops a coordinated bidding framework where multiple EVAs optimize their offers and exchange imbalanced liabilities through a Nash bargaining mechanism. In \cite{rt4}, a two-stage bidding approach is designed by integrating Stackelberg game-based optimization with chance constraints to coordinate dispatch decisions while addressing three-phase distribution imbalances.

While prior studies have contributed valuable insights into EVA participation in real-time electricity markets, three critical limitations remain, summarized as follows:
% Existing studies on EVAs have laid a strong foundation for EVA participation in the real-time market.

First, most existing EVA bidding strategies fail to reflect price-dependent operational preferences, limiting their responsiveness to market dynamics. Most existing works adopt a self-scheduling framework, wherein bids include only fixed energy quantities, indicating the amount of energy an EVA intends to purchase or supply, regardless of the actual market clearing price. Although this approach simplifies market participation, it prevents EVA from expressing its price preferences, often leading to suboptimal outcomes and limiting revenue potential \cite{xub}. To address these limitations, more flexible bidding mechanisms have been advocated, chief among which is economic bidding, where the EVA submits a set of price–quantity pairs that reflect its willingness to charge or discharge at different price levels \cite{ecobid}. This format enables EVA to respond adaptively to market prices and better align their operations with economic objectives. Fundamentally, each price–quantity pair in an economic bid reflects the EVA’s marginal valuation of energy at a given time, that is, the value (cost) of consuming (supplying) an additional unit of energy. However, unlike conventional generators whose marginal costs are typically derived from explicit and time-independent cost functions, an EVA’s marginal cost is shaped by implicit and inter-temporal opportunity costs. These arise from the potential gain or loss of future flexibility and are inherently difficult to model directly. Furthermore, they are also influenced by uncertain market prices and the dynamic availability of flexible EV resources, making them even harder to estimate accurately.

Second, the challenge of handling uncertainty, particularly in real-time electricity prices, remains inadequately addressed in existing EVA bidding strategies. A common approach is to rely on forecasts to inform bidding decisions. While forecast-based methods \cite{rt2,rt1,rt4} are simple and intuitive, their effectiveness heavily depends on both the accuracy of the forecasts and the prediction horizon. However, real-time prices are highly volatile and influenced by complex, non-stationary factors, causing prediction accuracy to degrade rapidly as the forecast horizon increases \cite{qiuj}. 
% which makes precise prediction notoriously difficult, 
To address these limitations, scenario-based methods \cite{dh1,dh2,dh3, rt3} aim to hedge against uncertainty by considering multiple future trajectories of uncertain parameters. 
% This approach offers probabilistic coverage of uncertainty, but faces a trade-off: a large number of scenarios increases modeling accuracy but leads to significant computational burden; conversely, reducing scenario size can compromise uncertainty coverage, potentially resulting in suboptimal bidding decisions. 
This approach offers probabilistic coverage of uncertainty but entails a trade-off between accuracy and efficiency: more scenarios increase computational cost, while fewer may lead to suboptimal outcomes.
% Robust optimization provides an alternative by seeking solutions that perform well under the worst-case realization of uncertainty. Although it is computationally efficient, its overly conservative nature often sacrifices economic performance. 
In addition, reinforcement learning has emerged as a data-driven approach that learns bidding policies through repeated interaction with the environment. Despite its flexibility, reinforcement learning remains limited by poor interpretability, difficulty in constraint enforcement, and high training costs, which may hinder its practical applicability to EVA bidding.

Third, the lack of dynamic, risk-averse bidding strategies remains a critical limitation in existing EVA real-time market participation frameworks. Given the inherent volatility of real-time prices, bidding decisions are inevitably exposed to substantial financial risk. However, many existing approaches adopt risk-neutral objectives, such as minimizing expected cost, which essentially assume that all potential outcomes are equally acceptable as long as their average is favorable. 
% This assumption may not hold in practice, as real-time electricity prices often exhibit heavy-tailed distributions.
This assumption may be inappropriate in real-time electricity markets, where price distributions are typically heavy-tailed in nature \cite{kim}. 
As a result, expectation-based objectives often fail to safeguard against rare but severe losses associated with extreme events (e.g., price spikes), leading to overly optimistic and financially vulnerable bidding strategies. To address this, several studies \cite{dh2,rt2} have incorporated risk measures into EVA bidding formulations. Among them, Conditional Value at Risk (CVaR) is widely used for capturing tail risk and balancing profitability with risk mitigation. Nevertheless, current CVaR-based bidding models suffer from several key limitations. Most notably, they usually rely on static risk measures, where risk is assessed once over a fixed time horizon based on pre-defined scenarios. This one-shot framework fails to reflect the sequential and evolving nature of real-time bidding, preventing risk re-evaluation as system states change. As a result, static formulations may lead to time-inconsistent strategies and degraded long-term performance. Moreover, they still rely on generating large scenario sets to approximate tail risks, which incurs significant computational overhead and complicates real-time deployment. 
% This assumption may not hold in practice, as real-time electricity prices often exhibit heavy-tailed distributions where rare but extreme price events can lead to substantial losses. As a result, strategies based solely on expected outcomes may underestimate financial exposure and fail to provide adequate protection against adverse scenarios. To mitigate this issue, recent studies have incorporated risk measures into EVA bidding models, among which Conditional Value at Risk (CVaR) is commonly employed to capture tail risk and balance profitability with robustness. Nevertheless, existing CVaR-based approaches are typically built on static formulations, where risk is evaluated once over a fixed planning horizon using pre-generated scenarios. This one-time assessment does not reflect the sequential nature of real-time bidding, where system states and uncertainties evolve continuously. The lack of dynamic risk assessment may result in time-inconsistent strategies that degrade long-term performance. Furthermore, the reliance on large scenario sets to approximate tail risks imposes considerable computational burdens, hindering practical implementation in real-time operational contexts.

To address the aforementioned research gaps, this paper develops an economic bidding strategy for EVA participating in real-time electricity markets. The goal is to minimize the total cost of the EVA under uncertainty in EV flexibility and real-time electricity prices, while satisfying operational constraints and complying with market bidding rules.
The main contributions of this work are summarized as follows:
% The proposed framework aims to minimize the total operational cost of the aggregator under price uncertainty and the evolving flexibility of EV resources, while satisfying operational constraints and complying with market bidding rules. 

{\begin{enumerate}
\item 
% We propose a marginal opportunity value–driven economic bidding strategy that enables EVA to express price-responsive charging and discharging preferences in real-time markets. Specifically, we formulate the EVA’s bidding problem as a Markov decision process (MDP) and use the associated value function to characterize the temporal and energy-dependent opportunity value. We derive a closed-form expression for the marginal opportunity cost from the derivative of the value function, and incorporate it into a price–quantity bidding structure. Additionally, we prove the convexity of the opportunity value function, ensuring that the resulting bidding strategy adheres to allowable bidding formats.
We propose a marginal opportunity value–driven economic bidding strategy that enables an EVA to express price-responsive charging and discharging preferences in real-time markets. By casting the sequential EVA scheduling problem as a Markov decision process (MDP), we explicitly link the inter-temporal opportunity value of energy to state-dependent value functions. On this basic, we derive the marginal opportunity values under both deterministic and stochastic price scenarios and embed them into the price-quantity bidding structure. This approach supports flexible integration of short-term forecast information with empirical knowledge from historical data, enhancing the strategy’s adaptability under varying levels of market information availability.
% We derive a closed-form expression for the marginal opportunity cost from the derivative of the value function, and incorporate it into a price–quantity bidding structure. 
% Additionally, we prove the convexity of the opportunity value function, ensuring that the resulting bidding strategy adheres to allowable bidding formats. 
% On this basis, we derive the marginal opportunity value from the derivative of the value function and incorporate it into the price-quantity bidding structure.
% Based on this, we further derive the probability distributions of marginal opportunity values for EVA under stochastic prices 
% In addition, it can flexibly incorporate short-term forecast information alongside empirical knowledge from historical data, further enhancing bidding performance in uncertain environments.
\item 
% We develop a computationally efficient, data-driven approach for real-time bidding under uncertainty. Instead of relying on extensive scenario sampling, the proposed method directly leverages the probabilistic structure of market uncertainty to train state- and time-dependent value functions within the MDP framework. Moreover, the proposed method is fully analytical, enabling direct computation of bidding decisions without the need for solver-based optimization. In addition, it can flexibly incorporate short-term forecast information alongside empirical knowledge from historical data, further enhancing bidding performance in uncertain environments.
We develop a fully analytical and computationally efficient approach for EVA bidding. Instead of relying on extensive scenario sampling, the proposed method directly leverages the probabilistic structure of electricity prices to construct the full probability distributions of marginal opportunity values across energy state. Accordingly, we further develop closed-form expressions for marginal charging values and discharging costs under both risk-neutral and risk-averse settings, enabling fast and interpretable bid generation through direct functional mapping, without relying on numerical solvers.
% enabling direct computation of bidding decisions without the need for solver-based optimization.
\item 
% We extend the conventional risk-neutral bidding framework into a dynamic, risk-averse formulation by integrating CVaR into the MDP setting. The model employs dynamic risk measures, where risk evaluation is recursively performed across decision stages using a time-consistent, nested formulation. This dynamic structure explicitly incorporates future risks into present decision-making, enabling adaptive and stage-wise management of the trade-off between expected profits and tail risks. Additionally, by analyzing the probability distribution and cumulative distribution of the value function, we derive a closed-form expression for the marginal opportunity cost under the dynamic risk measure, thus preserving both analytical tractability and computational efficiency.
% We introduce a dynamic risk-aware bidding formulation by embedding CVaR into the MDP framework. The model employs dynamic risk measures, where risk evaluation is recursively performed across decision stages using a time-consistent, nested formulation. This dynamic structure explicitly incorporates future risks into present decision-making, enabling adaptive and stage-wise management of the trade-off between expected profits and tail risks. 
We introduce a dynamic risk-aware bidding formulation by embedding CVaR into the MDP framework. Unlike static risk models that assess risk over a fixed horizon, our method employs dynamic risk measures \cite{drm1,drm2} to recursively evaluate risk across decision stages in a time-consistent, nested structure. This enables the EVA to incorporate evolving tail risks into current decisions and adaptively manage the trade-off between expected profits and downside protection.
% Additionally, by analyzing the probability distribution and cumulative distribution of the value function, we derive a closed-form expression for the marginal opportunity cost under the dynamic risk measure, thus preserving both analytical tractability and computational efficiency.
\end{enumerate}}
% \item
% We validate the proposed bidding strategy using real-world data on EV charging behavior (from Macao) and electricity prices (from NYISO). Simulation results demonstrate the proposed strategy’s superior economic performance, robustness under uncertainty, and real-time implementability.

This paper is organized as follows. Section \ref{Problem Statement} outlines the real-time market framework and introduces the modeling foundation for EVA bidding. Section \ref{ch3} develops a risk-neutral EVA bidding strategy grounded in marginal opportunity value analysis and Section \ref{ch4} extends this framework to the risk-averse case using dynamic risk measures. Section \ref{sec_case} presents the simulation results. Section \ref{sec_conclusion} concludes the paper.

\section{Problem Setting and Preliminaries}\label{Problem Statement}
\subsection{Market Overview}
We consider a real-time electricity market environment in which an EVA participates as a market entity. Through contractual agreements with EV owners, the EVA is authorized to centrally coordinate their charging activities. In this role, the EVA functions not only as a flexible electricity consumer but also as a potential energy supplier when V2G services are available. This dual-role prosumer model allows the EVA to strategically manipulate its power profile in response to market signals, leveraging flexibility to optimize economic outcomes while ensuring that individual charging requirements are met.

Real-time markets typically operate on short settlement intervals, such as every $5$ or $15$ minutes \cite{zhengshuo}. Participants are required to submit bids shortly before each interval, which are then cleared by the market operator based on system conditions. Although specific implementations vary across different markets, such as NYISO in the United States or NEM in Australia, the underlying bidding and settlement mechanisms remain largely consistent.
% wherein participating entities submit bids just ahead of each interval. These bids are then used by the system operator to determine dispatch schedules and market clearing prices in near real time. Although implementation details may vary across regions, such as the NYISO in the United States or the NEM in Australia, the underlying bidding and settlement mechanisms remain largely consistent.

\begin{figure}
	\centering
%	\vspace{-4mm}
	\includegraphics[width=0.95\columnwidth]{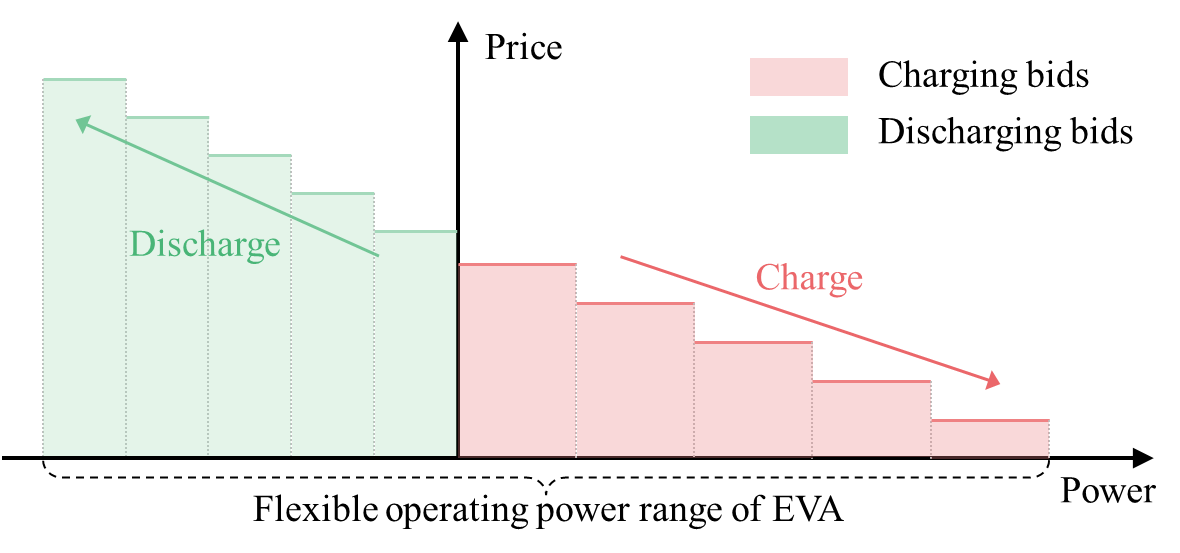}
	\vspace{-3mm}
    \caption{Illustrative example of EVA piecewise bidding.}
	\label{fig_bids}
\end{figure}

Within this framework, the EVA submits two-sided bids for each interval: one for charging (as a controllable load) and one for discharging (as a distributed generator). Fig. \ref{fig_bids} illustrates a representative bid structure for a single market interval. Each bid is structured as a stepwise curve composed of multiple price–quantity pairs, reflecting the EVA’s willingness to transact energy at different price levels. For charging bids, higher energy quantities are associated with lower bid prices, forming a monotonically decreasing demand curve. For discharging bids, the EVA offers to inject more energy only at higher prices, leading to a monotonically increasing supply curve.
% In this context, the EVA submits two-sided bids for every market interval: one reflecting its willingness to consume energy (charging bid) and the other representing its readiness to supply energy to the grid (discharging bid). Figure~1 illustrates the structure of a typical bid for a single market interval. Each bid is composed of a set of discrete price–quantity pairs that form a stepwise curve. The charging bid forms a monotonically decreasing demand curve, with larger quantities corresponding to lower bid prices—indicating that the EVA is willing to charge more only at lower prices. Conversely, the discharging bid represents a monotonically increasing supply curve, where higher quantities are offered only at higher prices, reflecting the opportunity cost of releasing stored energy.
% Charging bids are characterized by monotonically decreasing curves, where greater energy volumes are offered at lower prices. Conversely, discharging bids follow monotonically increasing patterns, with larger output volumes contingent on higher prices.

Once the market is cleared, the EVA must comply with its accepted bids. If the market price falls below a submitted charging price, the corresponding charging quantity is executed. If the price exceeds a discharging bid, the EVA discharges the associated amount. Otherwise, the EVA remains idle during the interval. In this study, we assume that the EVA operates as a price taker due to its relatively small market share, that is, its bids do not influence market clearing prices.
% In this study, we assume that the EVA operates as a price taker, meaning its bids do not influence the market clearing price—a reasonable assumption given its limited market size and participation volume.

\vspace{-1mm}
\subsection{Aggregated EV Flexibility Modeling}
The bidding behavior of an EVA fundamentally depends on its ability to quantify and manage the available flexibility of the EV fleet it controls. To participate effectively in the electricity market, the EVA must transform this distributed flexibility into tractable models that inform real-time bidding decisions.

\begin{figure}
	\centering
%	\vspace{-4mm}
	\includegraphics[width=0.95\columnwidth]{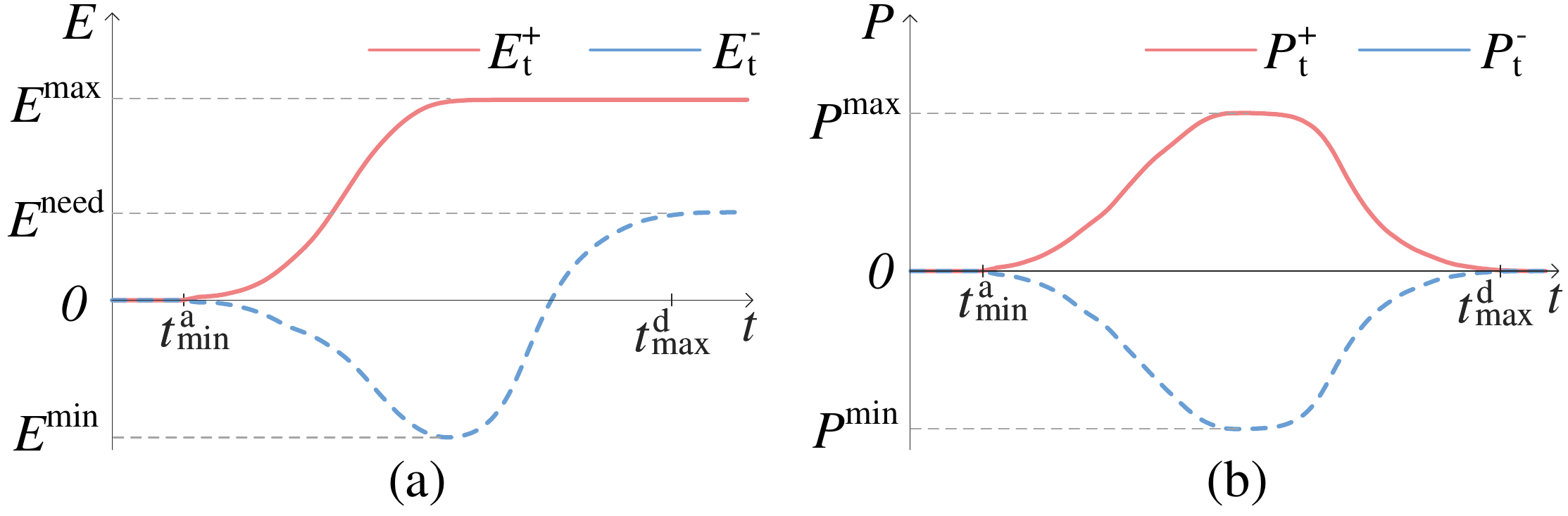}
	\vspace{-2mm}
    \caption{Aggregated energy and power boundaries of an EV fleet. In both subplots, $t_{\min }^{\text{a}}$ and $t_{\max }^{\text{d}}$ denote the earliest arrival and latest departure times of all EVs. In subplot (a), ${E^{{\text{max}}}}$, ${E^{{\text{min}}}}$ and ${E^{{\text{need}}}}$ denote the maximum, minimum, and required cumulative energy levels. In subplot (b), ${P^{{\text{max}}}}$ and ${P^{{\text{min}}}}$ are the maximum and minimum allowable charging/discharging power, determined by the physical limits of the EV fleet.}
	\label{fig_bound}
\end{figure}

Following our previous work \cite{zhang1}, we adopt a time-varying envelope-based model to characterize the feasible charging trajectories of the EV fleet. As illustrated in Fig. \ref{fig_bound}, the upper and lower energy envelopes, denoted by $E_t^ +$ and $E_t^ -$, represent the fastest and slowest admissible energy accumulation paths over time. These trajectories reflect the energy limits achievable when all EVs charge or defer charging at their respective maximum rates. 
Similarly, the EVA’s power capabilities are bounded by $P_t^ +$ and $P_t^ -$, corresponding to the fleet-level maximum charging and discharging power at each time slot. These envelopes are obtained by summing the individual energy and power limits of all connected EVs. Accordingly, the EVA’s operation must satisfy the following constraints:
\begin{align}
& 0 \leq p_t^{\text{c}} \leq P_t^ + , \quad \forall t, \label{cons1}\\
& P_t^ -  \leq p_t^{\text{d}} \leq 0, \quad \forall t, \\
& p_t^{\text{c}} \cdot p_t^{\text{d}} = 0, \quad \forall t, \label{cons3} \\
& {e_t} = {e_{t - 1}} + {p_t^{\text{c}}\eta  + p_t^{\text{d}}/\eta }, \quad \forall t, \label{cons4} \\
& E_t^ -  \leq {e_t} \leq E_t^ +, \quad \forall t, \label{cons5}
\end{align}
where $p_t^{\text{c}}$ and $p_t^{\text{d}}$ denote the EVA’s net charging and discharging power at time $t$, subject to the corresponding bounds $P_t^ +$ and $P_t^ -$. Constraint (\ref{cons3}) enforces mutual exclusivity of charging and discharging. The energy dynamics are governed by (\ref{cons4}), where ${e_t}$ represents the EVA’s energy level at time $t$, and $\eta$ is the charging/discharging efficiency. Constraint (\ref{cons5}) ensures that the EVA’s energy remains within feasible range throughout the operation horizon.
% This modeling framework captures the time-coupled flexibility of the EV fleet in a compact and analytically tractable form, making it a scalable foundation for integration into real-time market bidding frameworks.
% This modeling framework compactly captures the time-coupled flexibility of the fleet and supports scalable integration into real-time bidding strategies.

\subsection{Uncertainty-Aware EV Aggregation and Forecasting} 
As mentioned earlier, the EVA must estimate its available flexibility prior to each market interval. However, at the time of bidding, not all participating EVs are connected to the grid. To account for this, we partition the fleet into two groups: currently connected EVs, whose operational states are fully observable, and soon-to-arrive EVs, whose flexibility will become available during the upcoming market interval. For instance, when bidding at time $t$ for the interval $[t+1,t+2)$, only EVs already connected at $t$ contribute to the deterministic portion, while those expected to arrive within $[t,t+1)$ form the uncertain portion. EVs arriving after $t+1$ are excluded, as they do not impact the current decision window. 

The total flexibility used for bidding is then captured by combining these two components. In particular, the deterministic portion is updated from real-time data, while the uncertain portion is estimated using forecasts. Notably, compared to individual EV behavior, which is often irregular and stochastic, the aggregate dynamics of a large EV population tend to exhibit smooth and stable patterns, as shown in Fig. \ref{fig_bound}. This makes the aggregated energy and power boundaries well suited for short-term prediction using standard methods such as autoregressive models or exponential smoothing.

In addition, to further enhance robustness, chance-constrained formulations can be employed to account for prediction errors. These methods allow the EVA to conservatively adjust its operational boundaries, ensuring constraint satisfaction with high probability. Prior studies \cite{cc1, cc2}, including our earlier work \cite{zhang2}, have extensively explored this modeling approach and demonstrated its effectiveness in managing EVA's flexibility under uncertainty. For brevity, we omit the detailed formulation here and continue to use $ \{ E_t^ + ,E_t^ - ,P_t^ + ,P_t^ - \} $ to denote the updated, uncertainty-aware operational bounds throughout the remainder of this paper.

\section{Risk-Neutral EVA Bidding Strategy}\label{ch3}
\subsection{Problem Formulation and MDP Reformulation}
To design a market-oriented bidding strategy, we begin by analyzing the EVA’s underlying real-time power scheduling problem. Operating over a finite time horizon $t = 1,2, \dots ,T$, the EVA sequentially determines its aggregate charging and discharging power, $p_t^{\text{c}}$ and $p_t^{\text{d}}$, subject to the operational constraints defined in equations (\ref{cons1})–(\ref{cons5}). The objective is to minimize the expected cumulative cost across the horizon:

\vspace{-6mm}
\begin{align} 
\min \ \mathbb{E}\left[ c_t + \mathbb{E}[c_{t+1} + \dots + \mathbb{E}[c_T]]\right], \label{cons6}
\end{align}
where the immediate cost at time $t$, denoted ${c_t}$, is given by:
\begin{align} 
c_t = \pi_t^{\text{elec}} (p_t^{\text{c}} + p_t^{\text{d}}) - \pi^{\text{deg}} p_t^{\text{d}}. 
\end{align}
The first term reflects net expenditure or revenue from real-time market transactions, where $\pi_t^{\text{elec}}$ denotes the market clearing price. The second term accounts for EV battery degradation, modeled as proportional to the discharged energy, with a constant marginal degradation rate ${\pi ^{{\text{deg}}}}$.

To incorporate the sequential and uncertain nature of this problem, we reformulate it as a finite-horizon MDP. 
% To systematically handle temporal coupling and price uncertainty, we reformulate the scheduling problem as a finite-horizon MDP. 
% To explicitly capture the sequential, state-dependent nature of this problem under uncertainty, we reformulate it as a finite-horizon Markov Decision Process (MDP). 
At each decision epoch $t \in \mathcal{T} = \left\{ {1,2, \dots ,T} \right\}$, the system state is represented as:
\begin{align}
s_t = \left\{ \pi_t^{\text{elec}}, e_t, E_t^+, E_t^-, P_t^+, P_t^- \right\},
\end{align}
and the action is defined by:
\begin{align}
a_t = \left\{ p_t^{\text{c}}, p_t^{\text{d}} \right\}, \quad a_t \in \Omega_t,
\end{align}
where $\Omega_t$ denotes the feasible decision space derived from operational constraints (\ref{cons1})–(\ref{cons5}).

The stage reward is the negative of the immediate cost:
\begin{align}
r_t(s_t, a_t) = -\left[ \pi_t^{\text{elec}}(p_t^{\text{c}} + p_t^{\text{d}}) - \pi^{\text{deg}} p_t^{\text{d}} \right].
\end{align}

The system evolves according to two transition dynamics: a deterministic update of energy level based on charging/discharging actions, and a stochastic realization of future electricity price. 

Given this structure, the scheduling problem can then be solved recursively via the Bellman’s equation:
\begin{align}
{V_t}({e_t}) = \mathop {{\text{max}}}\limits_{{a_t} \in \Omega } \left\{ {{r_t}({s_t},{a_t}) + {\mathbb{E}_{\pi _{t + 1}^{{\text{elec}}}}}\left[ {{V_{t + 1}}({e_{t + 1}})} \right]} \right\}, \label{cons11}
\end{align}
% where ${V_t}({e_{t}})$ is the value function representing the maximum expected cumulative reward from time $t$ onward, given the current energy state ${e_{t}}$. 
% where ${V_t}({e_{t}})$ is the value function representing the maximum cumulative reward from time $t$ onward, given the current energy state ${e_{t}}$. 
where ${V_t}({e_{t}})$ is the value function, representing the opportunity value given the current energy state ${e_{t}}$, i.e., it reflects the economic benefit of deferring energy usage to exploit uncertain but potentially favorable future prices. Therefore, this recursive formulation captures the trade-off between immediate economic gains and the preservation of flexibility for future opportunities.
% This recursive formulation captures the trade-off between immediate economic gains and the deferred value of preserving flexibility for future use under uncertainty. 
To isolate the effect of future uncertainty, we introduce the post-decision value function \cite{powell} as:
% This recursive structure captures the trade-off between immediate economic gains and the strategic preservation of flexibility for future opportunities.
% This formulation captures the core trade-off faced by the EVA: balancing immediate market gains against the potential value of energy stored for future use.
\begin{align}
V_t^{\text{post}}(e_{t+1}) = \mathbb{E}_{\pi_{t+1}^{\text{elec}}} [ V_{t+1}(e_{t+1}) ],\label{cons12}
\end{align}
which evaluates the expected opportunity value after executing the current action and  transitioning to energy level ${e_{t+1}}$, but before the next market outcome is realized. 
% This formulation directly quantifies the opportunity value of stored energy, i.e., it reflects the economic benefit of deferring energy usage to exploit uncertain but potentially favorable future prices.
% In other words, this function enables the EVA to explicitly quantify the opportunity value of stored energy. 
% but prior to observing the next market outcome.
%This formulation allows the EVA to evaluate the opportunity cost or value of storing energy under uncertainty.
% It provides a convenient bridge between current control actions and future system performance and plays a central role in the design of the bidding strategy.

\subsection{Marginal Valuation of Charging and Discharging Actions}
% Building upon the MDP framework, this section establishes a formal link between the EVA’s internal decision-making and its market-oriented bidding behavior. Specifically, we leverage the post-decision value function to derive marginal valuations for charging and discharging decisions, which in turn reflect the EVA’s willingness to pay or be compensated for energy exchange under different system states. 
Building upon the MDP framework, we now derive the EVA’s internal marginal valuations for energy exchange decisions. These marginal values indicate how the EVA values incremental changes in energy and serve as a foundation for constructing bidding strategies. We begin by stating a structural property of the post-decision value function introduced in the previous section:
\begin{proposition} 
The post-decision value function $V_t^{{\text{post}}}(e)$ is concave in the energy level $e$, $\forall t \in \mathcal{T}$. \label{prop1}
\end{proposition}
Proof: see Appendix A in the supplementary material \cite{supplydocument}. 

The concavity of $V_t^{{\text{post}}}(e)$ implies that its derivative with respect to energy, denoted as $v_t^{{\text{post}}}(e) = \frac{{\partial V_t^{{\text{post}}}(e)}}{{\partial e}}$, is well-defined and monotonically non-increasing. Economically, $v_t^{{\text{post}}}(e)$ captures the expected marginal opportunity value of stored energy, i.e., the incremental benefit of storing one additional unit of energy at energy state $e$. Based on this foundation, we derive explicit expressions for the EVA’s marginal valuations of charging and discharging decisions.

\subsubsection{Marginal Value of Charging} Charging increases the EVA’s energy by $\eta  \cdot p_t^{\text{c}}$. Applying the chain rule to the post-decision value function yields: 
\begin{align} 
\nonumber {B_t^{\text{c}}({e_{t + 1}})} & = \frac{{\partial V_t^{{\text{post}}}({e_{t + 1}})}}{{\partial p_t^{\text{c}}}} = \frac{{\partial V_t^{{\text{post}}}({e_{t + 1}})}}{{\partial {e_{t + 1}}}}\frac{{\partial ({e_t} + p_t^{\text{c}}\eta )}}{{\partial p_t^{\text{c}}}} \\
& = \eta v_t^{{\text{post}}}({e_{t + 1}}). \label{cons13} 
\end{align}
This function calculates the maximum acceptable market price, $B_t^{\text{c}}$, for charging to remain profitable, given the resulting energy state ${e_{t + 1}}$.

\subsubsection{Marginal Cost of Discharging}
Likewise, discharging decreases the EVA’s energy by $\frac{1}{\eta } \cdot p_t^{\text{d}}$. The marginal cost of discharging incorporates both EV battery degradation and the opportunity cost of reduced future flexibility:
\begin{align}
\nonumber {B_t^{\text{d}}({e_{t + 1}})} & = \frac{{\partial ( - {\pi ^{{\text{deg}}}}p_t^{\text{d}} - V_t^{{\text{post}}}({e_{t + 1}}))}}{{\partial ( - p_t^{\text{d}})}} \\
\nonumber & = {\pi ^{{\text{deg}}}} - \frac{{\partial V_t^{{\text{post}}}({e_{t + 1}})}}{{\partial {e_{t + 1}}}}\frac{{\partial ({e_t} + p_t^{\text{d}}/\eta )}}{{\partial ( - p_t^{\text{d}})}} \\
& = {\pi ^{{\text{deg}}}} + \frac{{v_t^{{\text{post}}}({e_{t + 1}})}}{\eta }. \label{cons14}
\end{align}
This function calculates the minimum price, $B_t^{\text{d}}$, at which discharging becomes economically viable, given the associated energy state $e_{t + 1}$. 

These equations (\ref{cons13}) and (\ref{cons14}) provide an economic interpretation of the EVA’s charging and discharging preferences, forming the core of the EVA’s price-responsive behavior.
% These expressions provide economically meaningful thresholds that determine when it is profitable to buy (charge) or sell (discharge) electricity, forming the core of the EVA’s price-responsive behavior.
% serving as the basis for constructing economically interpretable, price-responsive bidding strategies.

\subsection{Market-Compliant Bid Construction} \label{C3S3}
% Having derived explicit expressions for the EVA’s marginal valuations of charging and discharging, we now describe how to translate these internal economic signals into market-compliant bid structures, i.e., stepwise monotonic price–quantity pairs. To this end, we propose a four-step four-step bid construction procedure, as illustrated in Figure 3.
To interface with real-time electricity markets, the EVA must convert continuous marginal valuations into discrete, market-compliant bids, typically structured as monotonic stepwise price–quantity pairs. We propose a four-step bid construction procedure.
%, as visualized in Fig. \ref{fig_bidgen}.

\begin{figure}
	\centering
%	\vspace{-4mm}
	\includegraphics[width=0.95\columnwidth]{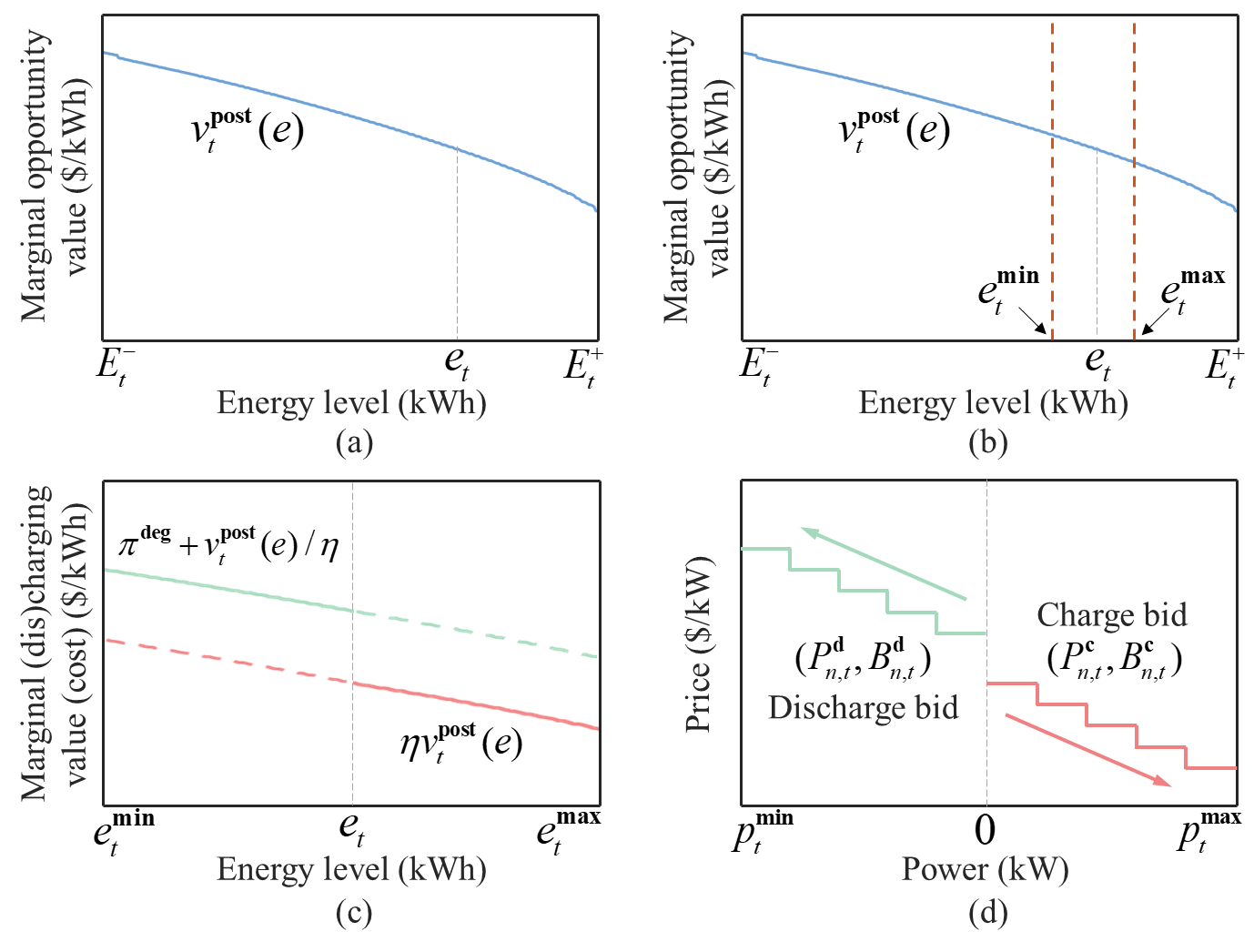}
	\vspace{-2mm}
    \caption{Illustrative example of EVA bids generation. (a) Post-decision marginal value function $v_t^{\text{post}}(e)$. (b) Feasible energy transition range $[e_t^{\min}, e_t^{\max}]$. (c) Marginal price functions for charging and discharging. 
    (d) Economic bid curves composed of price–quantity pairs.}
	\label{fig_bidgen}
 \vspace{-3mm}
\end{figure}

\subsubsection{Step1-Evaluate Expected Marginal Opportunity Values} To generate bids for the time period $t$, we begin by evaluating the post-decision marginal value function $v_t^{\text{post}}(e)$ across the admissible energy range $\left[ E_t^-, E_t^+ \right]$ (the detailed calculation method for $v_t^{\text{post}}(e)$ is deferred to subsection \ref{C3S4}). As illustrated in Fig. \ref{fig_bidgen}(a), this function is monotonically non-increasing, reflecting the diminishing marginal benefit of stored energy.

\subsubsection{Step2-Identify Feasible Energy Transitions} Given the initial energy level $e_t$, the range of reachable energy states during the bidding interval is constrained by the EVA’s aggregated power limits and energy boundaries:
\begin{align}
e_t^{{\text{max}}} = \min \left\{ {{e_t} + P_t^ + \eta ,E_t^ + } \right\},e_t^{{\text{min}}} = \max \left\{ {{e_t} + P_t^ - /\eta ,E_t^ - }\right\}. 
\end{align}
This range, depicted in Fig. \ref{fig_bidgen}(b), defines the energy flexibility that can be leveraged for market participation through either charging or discharging.

\subsubsection{Step3-Map to Marginal Price Functions} Over the intervals $[e_t, e_t^{\max}]$ for charging and $[e_t^{\min}, e_t]$ for discharging, we compute marginal valuations based on the previously derived equations (\ref{cons13})-(\ref{cons14}). As shown in Fig.~\ref{fig_bidgen}(c), the red solid line represents the marginal charging value and the green solid line represents the marginal discharging cost. The dashed lines represent marginal values in energy regions that lie outside the feasible direction of transition from $e_t$, and are thus excluded from bid construction.
% Over the intervals $[e_t, e_t^{\text{max}}]$ for charging and $[e_t^{\text{min}}, e_t]$ for discharging, we compute marginal valuations based on the previously derived equations (\ref{cons13})-(\ref{cons14}), as shown in Fig. \ref{fig_bidgen}(c).
% These continuous functions represent the EVA’s willingness to charging (pay) or discharging (accept) energy at different levels of state $e$ (Figure~3c).

\subsubsection{Step4-Generate Stepwise Bid Curves} 
To comply with market requirements, we discretize the feasible power range into $N$ equal segments and assign each segment a corresponding quantity and average marginal price:  
\begin{align}
& P_{n,t}^{\text{c}} = \frac{n}{N}p_t^{\max },\quad P_{n,t}^{\text{d}} = \frac{n}{N}p_t^{\min }, \\
& B_{n,t}^{\text{c}} = \eta  \cdot \frac{1}{{\Delta {e^{\text{c}}}}}\int_{{e_t} + (n - 1)\Delta {e^{\text{c}}}}^{{e_t} + n\Delta {e^{\text{c}}}} {v_t^{{\text{post}}}} (e), \\
& B_{n,t}^{\text{d}} = {\pi ^{{\text{deg}}}} + \frac{1}{\eta } \cdot \frac{1}{{\Delta {e^{\text{d}}}}}\int_{{e_t} - n\Delta {e^{\text{d}}}}^{{e_t} - (n - 1)\Delta {e^{\text{d}}}} {v_t^{{\text{post}}}} (e),  \\
& p_t^{{\text{max}}} = (e_t^{{\text{max}}} - {e_t})/\eta , \quad p_t^{{\text{min}}} = (e_t^{{\text{min}}} - {e_t})\eta , \\
& \Delta e^{\text{c}} = \frac{e_t^{\max} - e_t}{N}, \quad \Delta e^{\text{d}} = \frac{e_t - e_t^{\min}}{N}.
\end{align}
As depicted in Fig. \ref{fig_bidgen}(d), the constructed bids satisfy market monotonicity constraints:
\begin{align}
\text{Charging bids:} \quad B_1^{\text{c}} > B_2^{\text{c}} > ... > B_N^{\text{c}}, \\
\text{Discharging bids:} \quad B_1^{\text{d}} < B_2^{\text{d}} < ... < B_N^{\text{d}}.
\end{align}

%These bid curves reflect both physical constraints and economic priorities, enabling the EVA to participate in electricity markets in a strategically optimal manner.
These bid curves constitute the EVA’s final submission to the market and reflect both physical constraints and economic priorities, enabling the EVA to participate in electricity markets in a strategically optimal manner.
% encoding both operational feasibility and economic rationality in response to prevailing market prices.
% translating internal dynamic valuations into price signals that respect operational feasibility while maximizing economic responsiveness.

Upon market closure, the system operator clears all market participants' bids and sends the clearing price and quantity signals back to the EVA. During the dispatch period, the EVA fulfills its market commitments by properly distributing the total cleared power across the individual EVs. To achieve this, we adopt a heuristic power allocation strategy proposed in our previous work \cite{zhang1}, which allocates power based on the urgency of the EVs' charging demand and their remaining energy. Interested readers can refer to \cite{Lyu} for alternative power allocation algorithms. Finally, the EVA tracks the energy levels of each vehicle and updates the aggregation model's parameters, preparing for the next bidding cycle.

\subsection{Training of Post-Decision Marginal Value Function} \label{C3S4}

As established in the previous section, the post-decision marginal value function $v_t^{\text{post}}(e)$, which quantifies the expected marginal opportunity value of stored energy, is essential for guiding the EVA’s bidding decisions in real-time electricity markets. In this section, we develop a fully analytical offline training method for constructing $v_t^{\text{post}}(e)$ recursively across the entire state space, incorporating both short-term forecasts and long-term statistical knowledge.

\begin{algorithm}[!t] \label{alg}
    \caption{Offline Training Procedure for $v_t^{\text{post}}(e)$} 
    \label{alg:offline-training}
    \begin{footnotesize}
    \Input{The operational parameters of EVA; historical electricity price distributions $f_t^{\text{elec}}(\cdot)$ and $F_t^{\text{elec}}(\cdot)$.}
    \Output{Well-trained post-decision marginal value functions $v_t^{\text{post}}(e)$ for all $t$ and $e$.}
    \textbf{Procedure:} \\
    \textbf{Initialization:} Set the forecast horizon $H$; discretize the EVA's energy level $e$ into $M$ segments uniformly; specify the terminal marginal value function $v_T^{\text{post}}(e)$. \\
    {\text{Backward recursion with historical distribution:}} \\
    \For{$t = T, \dots, H+1 $}{
        \For{$m = 1, \dots, M $}{
                Update $v_t^{\text{post}}(e_m)$ from $v_{t+1}^{\text{post}}(e_m)$ using Equation (\ref{cons25}); \\
                Apply boundary correction: \\
                $v_t^{\text{post}}(e_m) = 
                \begin{cases}
                \chi, & \text{if } e_m < E_t^-, \\
                0, & \text{if } e_m \geq E_t^+,
                \end{cases}$ where $\chi \gg 0$\;
            }
            }
    {\text{Backward recursion with short-term forecast:}} \\       
    \For{$t = H, \dots, 1$}{
        \For{$m = 1, \dots, M $}{
                Update $v_t^{\text{post}}(e_m)$ from $v_{t+1}^{\text{post}}(e_m)$ using Equation (\ref{cons23}); \\
                Apply boundary correction as above. \\
            }
            }        
    \end{footnotesize}
\end{algorithm}
% 3. **Boundary Corrections**: The piecewise structure of \( v_t^{\text{post}}(e) \) requires corrections at the boundaries of each interval. For example, if \( \pi_t^{\text{elec}} \) exceeds certain thresholds, the corresponding energy level will be adjusted according to the given constraints, ensuring that the solution respects the energy bounds and pricing structure.

\subsubsection{Recursive Updates Under Deterministic Price Forecasts} 
% When future electricity prices are deterministically known, the marginal opportunity value function becomes deterministic as well, i.e., $v_t^{{\text{post}}}({e_{t + 1}}) = {v_{t + 1}}({e_{t + 1}})$. While this assumption may not realistic in practice, it provides useful insight and forms the theoretical foundation for later extensions.
When future electricity prices are deterministically known, the marginal opportunity value becomes deterministic as well. Specifically, according to equation~(\ref{cons12}), the post-decision marginal value function now coincides with the marginal value function, that is, $v_t^{\text{post}}(e_{t+1}) = v_{t+1}(e_{t+1})$. Although the assumption of deterministic prices may not hold in practice, it provides useful theoretical insight and serves as a basis for later extensions.
% When deterministic short-term forecasts of electricity prices are available, they can provide valuable guidance for shaping the EVA’s real-time bidding behavior. 
% In this case, we assume that future prices within a forecast window are known deterministically. 
% This enables the post-decision marginal value function to be computed recursively in closed form based on the anticipated price path.
\begin{proposition} \label{prop2}
Given deterministic electricity prices within a forecast window, the post-decision marginal value function satisfies the following backward recursion: 
\begin{align} 
v_t^{\text{post}}(e) = G^{\text{d}}[v_{t+1}^{\text{post}}(e)]. \label{cons23}
\end{align} 
\end{proposition}
The detailed formulation of the operator $G^{\text{d}}[\cdot]$ and the proof are provided in Appendix B of the supplementary material \cite{supplydocument}.

\subsubsection{Recursive Updates Using Historical Price Distributions} 
In practice, perfect foresight is rarely attainable. While short-term forecasts may capture some trends, their accuracy typically degrades rapidly over extended horizons. When electricity prices are uncertain, the marginal opportunity value also becomes stochastic. In this study, we model the real-time electricity prices $\pi_t^{\text{elec}}$ as stage-wise independent stochastic variables following probability distribution $f_t^{\text{elec}}(\cdot)$ and cumulative distribution $F_t^{\text{elec}}(\cdot)$, both derived from historical market data. Based on this setup, we derive the probability distribution of the EVA’s marginal opportunity value, i.e., the derivative of the value function $V_t(e)$ with respect to energy, denoted by $v_t(e)$. The distribution is given by:
\begin{align}
\notag & \Pr[v_t(e) = x] = \\
& \begin{cases}
F_t^{\text{elec}}[x\eta],  \quad  \text{if } x = v_t^{\text{post}}(e + \eta P_t^+) \\
\eta f_t^{\text{elec}}[x\eta], \quad \text{if } v_t^{\text{post}}(e + \eta P_t^+) < x < v_t^{\text{post}}(e) \\
F_t^{\text{elec}}\left[\frac{x}{\eta} + \pi^{\text{deg}}\right] - F_t^{\text{elec}}[x\eta],  \quad\text{if } x = v_t^{\text{post}}(e) \\
\frac{1}{\eta} f_t^{\text{elec}}\!\left[\frac{x}{\eta} \! + \!\pi^{\text{deg}}\right], ~~ \text{if } v_t^{\text{post}}(e) \! < x < \!v_t^{\text{post}}(e + \! \frac{P_t^-}{\eta}) \\
1 - F_t^{\text{elec}}\left[\frac{x}{\eta} + \pi^{\text{deg}}\right], \quad \text{if } x = v_t^{\text{post}}(e + \frac{P_t^-}{\eta}) \\
0,  \quad \text{else}
\end{cases},
\end{align}

The proof is provided in Appendix C. Then, taking the expectation of $v_t(e)$ over this distribution yields the post-decision marginal value function.
\begin{proposition} 
Given $f_t^{{\text{elec}}}( \cdot )$ and $F_t^{{\text{elec}}}( \cdot )$ estimated from historical prices, the post-decision marginal value function admits the following recursive update:
\begin{align} 
v_t^{\text{post}}(e) = G^{\text{s}}[v_{t+1}^{\text{post}}(e)]. \label{cons25}
\end{align} 
\end{proposition}
% The proof is provided in Appendix C.
The explicit formulation of the operator $G^{\text{s}}[\cdot]$ and the proof are provided in Appendix C of the supplementary material \cite{supplydocument}.
% This proposition is derived by analytically establishing the probability distribution functions for the marginal value function ${v_t}(e)$, i.e., the marginal opportunity value of EVA, and then calculating the expectation integrals. Given the complexity and space constraints, the explicit formulation of the operator $G^{\text{s}}[\cdot]$ and a detailed proof are deferred to Appendix~C.

Based on the above formulations, we design a hybrid training framework that integrates both historical data and forecast information. This framework enables the EVA to adaptively leverage accurate predictions when available, while maintaining robustness under uncertainty through historical distributions. Starting from the predefined terminal condition $v_T^{\text{post}}(e)$ that encodes the final energy requirement (see Appendix A), the training process proceeds backward in time, recursively computing $v_t^{\text{post}}(e)$ across all energy states and time stages. The full training procedure is outlined in Algorithm~\ref{alg:offline-training}.

The proposed training method offers several advantages. First, it enables flexible integration of forecast and historical information through a tunable parameter $H$, which allows adaptation from fully historical data-driven ($H=0$) to fully forecast-driven ($H=T$) training. Second, the algorithm scales linearly with the time horizon $T$ and the number of energy discretization points $M$, resulting in a overall computational complexity of $\mathcal{O}(TM)$. Specifically, it performs $M$ updates of $v_t^{\text{post}}(e)$ at each time step, thereby effectively avoiding the curse of dimensionality. Lastly, the entire training is closed-form and solver-free, eliminating the need for numerical optimization and enabling fast, scalable deployment. 
% Third, memory usage remains constant with respect to $T$ since only the current value function array needs to be retained at each step. 

\section{Risk-Averse EVA Bidding Strategy}\label{ch4}
\subsection{Problem Formulation under Risk-Averse Case} 
In real-time electricity markets, an EVA is exposed to significant financial risks due to high price volatility. Traditional risk-neutral strategies, which focus solely on maximizing expected profits, fail to account for rare but severe events (such as price spikes) that could result in substantial financial losses. Therefore, it is essential to develop a more robust, risk-averse bidding strategy that safeguards the EVA against downside risks and ensures financial sustainability.

Following a similar analysis structure as in Section~\ref{ch3}, we first extend the risk-neutral objective in equation~(\ref{cons6}) to a risk-averse formulation:
\begin{align}
\max \ \rho_{t} \left[ r_t + \rho_{t+1} \left[ r_{t + 1} + \dots + \rho_{T} \left[ r_T \right] \right] \right],
\end{align}
where $\rho_{t}$ denotes the dynamic risk measure operator applied at each decision stage. In short, dynamic risk measures are an extension of traditional risk metrics, incorporating the evolving nature of risk over sequential decision-making processes. Interested readers can refer to \cite{drm1,drm2} for further details on dynamic risk measures. This formulation enables the EVA to dynamically adapt its decisions to emerging market conditions, rather than relying on a static, one-time evaluation.
% This formulation allows EVA to dynamically adjust its strategy as market conditions evolve, rather than relying on a static, one-time evaluation.

% One widely used dynamic risk measure is Conditional Value at Risk (CVaR), which is favored for its time-consistent properties. 
Among various risk measures, the Conditional Value-at-Risk (CVaR) is widely used due to its convexity, coherence, and time consistency properties \cite{ope}. CVaR is particularly useful in contexts where the risk of extreme market events (i.e., price spikes) needs to be controlled. Mathematically, CVaR at a given confidence level $\alpha$ is defined as:
% \begin{align}
% & \text{CVaR}_\alpha(X) \!\!= \!\mathbb{E}[X | X \! \leq \!\text{VaR}_\alpha] \!= \!\frac{1}{1-\alpha} \!\int_{-\infty}^{\text{VaR}_\alpha} \!x f(x) \, dx, \\
% & \text{VaR}_\alpha(X) = \max \{ Y | \Pr(X \leq Y) \leq 1 - \alpha \},
% \end{align}
\begin{align}
& \text{CVaR}_\alpha(X) \! = \! \mathbb{E}[X | X  \! \leq \! \text{VaR}_\alpha] \!=\! \frac{1}{1\!-\!\alpha} \! \int_{-\infty}^{\text{VaR}_\alpha} x f(x) \, dx, \\
& \text{VaR}_\alpha(X) = \max \{ Y | \Pr(X \leq Y) \leq 1 - \alpha \},
\end{align}
where $\text{VaR}_\alpha$ is the Value-at-Risk (VaR) at confidence level $\alpha$, defined as the maximum value such that the probability of profits being below this value is no more than $1 - \alpha$. CVaR then quantifies the expected profit conditional on outcomes that fall below the VaR threshold, emphasizing tail risk management. By incorporating CVaR, the EVA can effectively manage the risks associated with extreme market conditions, avoiding excessive financial losses.

To strike a balance between profitability and risk mitigation, a combination of the expected profit and CVaR is commonly employed \cite{math}, resulting in the risk measure: 
% we combine the expected profit and CVaR in a weighted manner, as is commonly done in many risk-averse models. The resulting risk measure $\rho_{t}$ is given by:
\begin{align}
\rho(X) = (1 - \lambda) \cdot \mathbb{E}[X] + \lambda \cdot \text{CVaR}_\alpha(X),
\end{align}
where $\lambda$ is a weighting parameter that controls the trade-off between expected profit and risk aversion. When $\lambda = 0$, EVA adopts a purely risk-neutral strategy, while $\lambda = 1$ corresponds to a fully risk-averse approach. This combined formulation enables the EVA to tailor its bidding strategy according to its risk preference.

\subsection{Risk-Averse MDP Formulation and Bid Construction}
Building on the risk-averse scheduling framework, we now extend the standard MDP formulation (\ref{cons11})–(\ref{cons12}) by incorporating dynamic risk measures. The resulting risk-averse MDP is given by:
\begin{align}
& V_t^{\text{R}}(e_t) = \max_{a_t \in \Omega} \left\{ r_t(s_t, a_t) + \rho \left[ V_{t+1}^{\text{R}}(e_{t+1}) \right] \right\}, \\
& V_t^{\text{post.R}}(e_{t+1}) = \rho \left[ V_{t+1}^{\text{R}}(e_{t+1}) \right],
\end{align}
where $V_t^{\text{R}}(e_t)$ and $V_t^{\text{post,R}}(e_{t+1})$ denote the risk-averse value function and post-decision value function, respectively. In particular, $V_t^{\text{post,R}}(e_{t+1})$ characterizes the risk-adjusted opportunity value of stored energy, accounting for both the expected future benefits and the potential risks.
% reflecting the expected gains of deferring energy usage while accounting for potential future risks.
Accordingly, we establish a key structural property:
\begin{proposition}
The risk-averse post-decision value function $V_t^{\text{post.R}}(e)$ is concave in the energy level $e$, $\forall t \in \mathcal{T}$.
\end{proposition}

The proof follows similar reasoning to Proposition \ref{prop1}. Briefly, since the immediate reward $r_t$ is linear in $e$, and the dynamic risk measure $\rho[\cdot]$ preserves concavity when applied to concave functions, the recursive application of Bellman’s equation maintains concavity over time. Consequently, the post-decision marginal value function under risk aversion, $v_t^{\text{post,R}}(e) = \frac{\partial V_t^{\text{post,R}}(e)}{\partial e}$ remains monotonically non-increasing.
% due to the linearity of the immediate reward $r_t$ with respect to the energy state $e$, and given that the dynamic risk measure $\rho[\cdot]$ preserves concavity when applied to concave functions, the recursive structure defined by Bellman’s equation ensures the preservation of concavity across all time stages.

Based on this property, the EVA's risk-averse bidding strategy can be generated by adapting a similar procedure developed in Section \ref{C3S3}. In short, the EVA first determines marginal opportunity values under risk aversion, and subsequently transforms these valuations into market-compliant, monotonic stepwise price–quantity bids.

The central remaining challenge is the computation of the risk-averse post-decision marginal value function $v_t^{\text{post,R}}(e)$, which governs the EVA’s valuation of energy increments under dynamic risk considerations. To this end, we develop an efficient analytical method for computing $v_t^{\text{post,R}}(e)$.
% which fundamentally guides the EVA’s bid formation under risk considerations. 

Specifically, under deterministic price forecast scenarios, risk adjustments become redundant due to the absence of uncertainty and thus no additional risk premium is necessary. Therefore, the training of $v_t^{\text{post,R}}(e)$ coincides exactly with the risk-neutral case described in Proposition \ref{prop2}. While under price uncertainty, we establish the following proposition, extending the recursive computation to the risk-averse context:
\begin{proposition} 
Given $f_t^{\text{elec}}(\cdot)$ and $F_t^{\text{elec}}(\cdot)$ derived from historical prices, the risk-averse post-decision marginal value function admits the following recursive update:
\begin{align}
v_t^{{\text{post}}{\text{.R}}}(e) = {G^{\text{R}}}[v_{t + 1}^{{\text{post}}{\text{.R}}}(e)].
\end{align}
\end{proposition}
The explicit form of the operator ${G^{\text{R}}}[\cdot]$ and a detailed proof are provided in Appendix~D of the supplementary material \cite{supplydocument}.
% This proposition is derived by analytically establishing the probability density and cumulative distribution functions for the risk-averse marginal value function $v_t^{\text{R}}(e)$, and then calculating the expectation and CVaR integrals. Given the complexity and space constraints, the explicit form of the operator ${G^{\text{R}}}[\cdot]$ and a detailed proof are deferred to Appendix~D.
Leveraging the recursive updates described above, the risk-averse post-decision marginal value function $v_t^{\text{post,R}}(e)$ can be efficiently computed for all energy states and decision stages. This can be achieved through a procedure analogous to Algorithm \ref{alg}, thus providing EVA with the necessary information to formulate dynamic, risk-averse bidding decisions in real-time electricity markets.
% Specifically, the EVA first proceeds backward from the terminal stage using the operator $G^{\text{R}}[\cdot]$ when historical distributions are used, and switches to $G^{\text{d}}[\cdot]$ when forecast information becomes available within a predefined forecast horizon $H$.
% Specifically, beginning from an appropriately defined terminal marginal condition, the EVA recursively applies the operators $G^{\text{d}}[\cdot]$ and $G^{\text{R}}[\cdot]$ backward in time, transitioning smoothly from forecast-based deterministic updates (if available) to probabilistic updates based on historical distributions.

\section{Case Studies}\label{sec_case}

\subsection{Case Settings}\label{sec_case_para}
\begin{figure}
\vspace{-2mm}
\centering
\includegraphics[width=1\columnwidth]{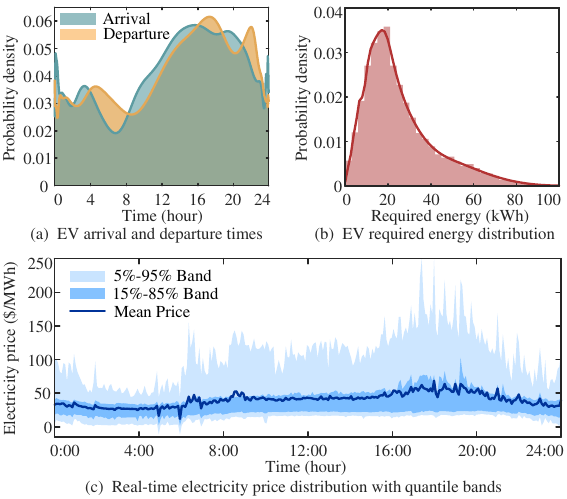}
	\vspace{-6mm}
	\caption{The statistical characteristics of EV charging and electricity prices.}
    \label{fig_characteristics}
\end{figure}

To verify the performance of the proposed bidding strategies, we conduct a simulation study that integrates real-world EV charging behavior with real-time electricity market conditions. The EV charging dataset originates from public charging stations in Macao SAR, China during March $2021$. After data preprocessing, a total of $10,040$ valid charging records are retained, each containing the EV’s plug-in time, plug-out time, and required energy to be charged. 
The distributions of these variables are shown in Fig. \ref{fig_characteristics}(a)–(b).
%Each record includes information such as EV’s plug-in time, plug-out time, and energy charged. 
% The distributions of arrival/departure times and required energy are shown in Fig.~4(a)–(b).
%The distributions of arrival and departure times, as well as required energy, are illustrated in Fig.~4(a)-(b). 
For EVA modeling purposes, we adopt the Tesla Model~S as the representative EV, with a battery capacity of $100$~kWh. Each charging pile is assumed to have a rated power of $50$~kW and supports bidirectional operation.
% The rated power of each charging pile is set to 50~kW and assumed to support V2G operation. 
The electricity price data is sourced from the NYISO $2018$ real-time market \cite{nyiso}. Fig. \ref{fig_characteristics}(c) illustrates the intra-day variation in real-time electricity prices. In our implementation, electricity prices are modeled as stage-wise independent Gaussian random variables. Eleven months of historical data are used to estimate the time-varying mean and variance, which are then used to construct the probability distribution function $f_t^{\text{elec}}(\cdot)$ and cumulative distribution function $F_t^{\text{elec}}(\cdot)$ for training the marginal value functions under both risk-neutral and risk-averse settings. The remaining data is used for validation.

The simulation spans $31$ consecutive days with a time resolution of five minutes, consistent with the NYISO real-time market settlement interval. 
%Each time step includes bid generation, market clearing, and intra-EVA power allocation. 
Cross-day charging behavior is supported, allowing EVs to plug in and out across different calendar days. System parameters are configured as follows: charging and discharging efficiency $\eta = 0.93$; marginal battery degradation cost $\pi^{\text{deg}} = 10~\$/\text{MWh}$; number of energy discretization levels $M = 1000$; number of bid steps $N = 5$.
% Both bidding and dispatching decisions are made accordingly in each interval.
% The charging and discharging efficiency is set to $\eta = 0.93$, and the marginal battery degradation cost is assumed to be $\pi^{\text{deg}} = 10~\$/\text{MWh}$. The energy level is discretized into $M = 1000$ segments to support value function approximation, and the bidding strategy is implemented using $N = 5$ discrete price–quantity pairs. 

\subsection{Performance Evaluation under Risk-Neutral Settings} 
This part evaluates the performance of the proposed EVA bidding strategies under risk-neutral settings. To benchmark their performance, we consider three comparative scenarios as follows:

\textbf{B1:} Uncoordinated charging. EVs begin charging at full rated power upon arrival and stop once their energy demand is fulfilled. This setup represents typical unoptimized charging behavior that ignores arbitrage potential.

\textbf{B2:} Ideal case. Each EV individually solves a deterministic optimization problem assuming perfect knowledge of future electricity prices. While infeasible in practice, this scenario serves as a theoretical upper bound on achievable economic efficiency.
% Although infeasible in practice, this scenario provides a theoretical baseline for maximum achievable economic efficiency under full information.

\textbf{B3:} Day-ahead price-informed scheduling. Each EV optimizes its charging schedule using day-ahead electricity prices, which are publicly available in advance. This represents a practical yet inherently suboptimal decision-making strategy.

It is worth emphasizing that all three benchmarks (B1–B3) operate outside the bidding-based market mechanism and assume direct control over individual EV charging profiles. Additionally, we evaluate the following two strategies apply the proposed EVA bidding framework:

\textbf{S1:} The proposed EVA bidding strategy, which constructs price–quantity bids based on the historical probability distributions of real-time electricity prices.

\textbf{S2:} A variant of S1, which models real-time prices through a combination of known day-ahead prices and historical probability distributions of deviations between real-time and day-ahead prices.

\begin{figure}
\centering
\includegraphics[width=0.75\columnwidth]{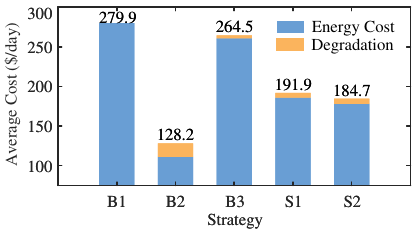}
	\vspace{-2mm}
	\caption{Operational cost comparison under different strategies.}
    \label{fig_cost}
\end{figure}

% Fig.~5 presents the total operational costs across all strategies. As expected, the uncoordinated benchmark B1 exhibits the highest operational cost at \$$8678$, due to its inability to exploit favorable price conditions. In contrast, the ideal case B2 achieves a total cost of \$$3973$, representing a $54.2$\% reduction relative to B1. This significant cost saving arises from fully leveraging the temporal flexibility of EVs to shift energy consumption to low-price periods. Notably, B2 also incurs the highest battery degradation cost, due to frequent and aggressive battery cycling aimed at capturing nearly all available arbitrage opportunities. In addition, benchmark B3 performs only marginally better than B1, achieving a cost reduction of about $5.5$\%. This reveals the limitations of using day-ahead prices as a proxy for real-time price dynamics, as it fails to capture intraday price volatility, leading to suboptimal or even misguided scheduling decisions. 
Fig.~\ref{fig_cost} presents the average daily operational costs across all strategies. As expected, the uncoordinated benchmark B1 exhibits the highest operational cost at \$$279.93$ per day, due to its inability to exploit favorable price conditions. In contrast, the ideal case B2 achieves a daily cost of \$$128.17$, representing a $54.2$\% reduction relative to B1. This significant cost saving arises from fully leveraging the temporal flexibility of EVs to shift energy consumption to low-price periods. Notably, B2 also incurs the highest battery degradation cost, due to frequent and aggressive battery cycling aimed at capturing nearly all available arbitrage opportunities. In addition, benchmark B3 performs only marginally better than B1, achieving a cost reduction of about $5.5$\%. This reveals the limitations of using day-ahead prices as a proxy for real-time price dynamics, as it fails to capture intraday price volatility, leading to suboptimal or even misguided scheduling decisions. 

% In comparison, both proposed EVA bidding strategies yield substantial performance gains. In particular, S1 achieves a $31.5$\% cost reduction relative to B1, lowering the operational cost to \$$5948$. This result demonstrates the effectiveness of using historical real-time price distributions to inform bid construction under uncertainty. In addition, S2 further reduces the cost to \$$5724$, achieving an additional $3.8$\% improvement over S1. This gain is attributed to anchoring the uncertain real-time price distribution around observed day-ahead prices, thereby narrowing the uncertainty range and enhancing the accuracy of marginal value estimation.
In comparison, both proposed EVA bidding strategies yield substantial performance gains. In particular, S1 achieves a $31.5$\% cost reduction relative to B1, lowering the daily cost to \$$191.87$. This result demonstrates the effectiveness of using historical real-time price distributions to inform bid construction under uncertainty. In addition, S2 further reduces the cost to \$$184.67$, achieving an additional $3.8$\% improvement over S1. This gain is attributed to anchoring the uncertain real-time price distribution around observed day-ahead prices, thereby narrowing the uncertainty range and enhancing the accuracy of marginal value estimation.

\subsection{Value of Short-Term Forecasts in EVA Bidding} 
This part investigates how the availability of short-term forecast information influences the performance of the proposed EVA bidding strategy. Building on the hybrid training framework described in Section \ref{C3S4}, we consider an extended version of strategy S2. In this variant, real-time electricity prices are assumed to be perfectly known within a specified forecast window, while price uncertainty beyond this window is still modeled using historical distributions of price deviations. Although perfect forecast assumption is optimistic, it enables a clear quantification of their theoretical contribution to performance enhancement.

\begin{figure}
\centering
\includegraphics[width=0.7\columnwidth]{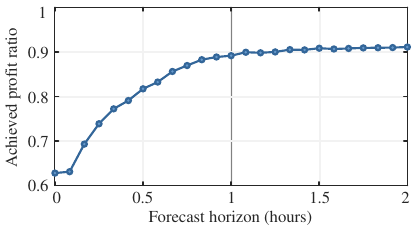}
	\vspace{-2mm}
	\caption{Profit ratio under different forecast horizons.}
    \label{fig_mpc}
\end{figure}

We define a normalized metric referred to as the profit ratio, which measures the proportion of total potential cost savings (i.e., the difference between the costs under the uncoordinated benchmark B1 and the ideal benchmark B2) captured by the strategy. The index ranges from $0$ to $1$, with higher values indicating closer alignment with the ideal cost outcome.

Simulation results are presented in Fig. \ref{fig_mpc}. As the forecast horizon increases, the profit ratio rises steadily, reflecting improved bidding performance. In particular, with a forecast window of only $30$ minutes, over $80$\% of the total potential cost savings are realized. As the forecast horizon extends to one hour, the realized profit ratio approaches $90$\%, and additional forecast length provides only marginal benefits. These results demonstrate that even limited forecast information can significantly enhance the performance of EVA bidding when incorporated into the proposed framework. In addition, a small residual gap remains relative to the ideal benchmark. This is primarily due to structural approximations in the aggregate EVA model and the discretized nature of the bid construction. Moreover, the realized flexibility of the EVA is influenced by the power allocation scheme that governs how aggregate instructions are distributed among individual EVs, which also contributes to the performance gap \cite{Lyu}.

\subsection{Comparative Analysis of Risk-Neutral and Risk-Averse Bidding Strategies} 
This part analyzes how incorporating risk aversion affects EVA bidding behavior and economic outcomes. The risk-averse strategy S3 extends the formulation of S2 by integrating the CVaR-based dynamic risk measure introduced in Section \ref{ch4}. Table \ref{tab} presents the operational costs under various degrees of risk aversion. Compared to the risk-neutral baseline S2, all tested risk-averse scenarios exhibit higher total costs but notably reduced battery degradation costs. This trend intensifies as either the risk weighting parameter $\lambda$ or the CVaR confidence level $\alpha$ increases, reflecting more conservative operational decisions.

\begin{table}[]
\footnotesize
\caption{One-Month simulation results under different risk settings}
\label{tab}
\begin{tabular}{c|cc|ccc|c}
\hline\hline
                    & $\lambda$ & $\alpha$ & \tabincell{c}{Energy \\ cost (\$) } & \tabincell{c}{Degradation \\ cost (\$) } & \tabincell{c}{Total \\ cost (\$) } & \tabincell{c}{CPU  time \\ ({\text{mins}}) }  \\ \hline
S2                  & 0                      & -                     & 5517.02         & 207.70               & 5724.71         & $\approx 7$                   \\
\multirow{4}{*}{S3} & 0.1                    & 95\%                  & 5636.04         & 145.00               & 5781.04         & \multirow{4}{*}{$\approx 28$} \\
                    & 0.1                    & 99\%                  & 5652.89         & 139.58               & 5792.47         &                       \\
                    & 0.2                    & 95\%                  & 5821.15         & 106.92               & 5928.07         &                       \\
                    & 0.2                    & 99\%                  & 5840.88         & 99.44                & 5940.32         &                       \\ \hline\hline
\end{tabular}
\end{table}

\begin{figure}
	\centering
	\includegraphics[width=1\columnwidth]{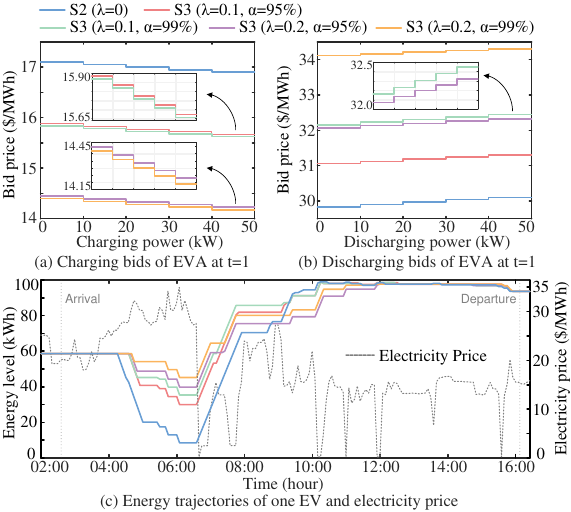}
	\vspace{-5mm}
	\caption{Impact of risk aversion on EVA bidding behavior and charging trajectories.}
    \label{fig_soc}
\end{figure}

To illustrate the underlying mechanism behind these cost differences, Figs. \ref{fig_soc}(a)-(b) compares the EVA bidding curves at the initial market interval ($t=1$) for different risk preferences. A clear pattern emerges: as EVA adopts higher risk aversion (larger $\lambda$ or $\alpha$), the bidding prices for charging decrease, while those for discharging increase. Intuitively, under risk-neutral conditions, bids represent the expected marginal opportunity value of energy increments, which might inadequately capture adverse market scenarios. In contrast, risk-averse bids explicitly incorporate tail risks by combining the expected marginal opportunity values with their CVaR. Consequently, charging bids are lowered to avoid potentially overvaluing future energy storage opportunities, while discharging bids increase to guard against the risk of realizing severe economic losses under unfavorable (though possibly rare) market conditions.
% charging bids decrease to guard against scenarios where the stored energy fails to yield sufficient future value, 
Moreover, it is notable that charging bids exhibit relatively minor differences between confidence levels $\alpha=95\%$ and $\alpha=99\%$. This suggests that the lower-tail distributions (i.e., worst 1\% and 5\% of marginal valuations) may exhibit limited differences at the examined market interval, resulting in similar conservative bids. 
% In contrast, the discharging bids in Fig.~7(b) display more evident variations across different risk parameters. This is because, for discharging decisions, EVA evaluates the upper-tail portion of marginal opportunity costs (i.e., the highest 1% and 5%), and the difference between these extreme values can be more pronounced, leading to clearer distinctions in discharging bids under varying levels of risk aversion.
% a heavy concentration of the lower-tail distribution at the examined market interval, where the 1% and 5% conditional expectations are numerically close.

To provide additional operational insights, we analyze the charging trajectory of a representative EV arriving at $2:40$ with an initial battery level of $58.7$ kWh, aiming to reach a target of $93.7$ kWh by $16:15$. Fig. \ref{fig_soc}(c) compares the corresponding charging paths across strategies with different risk preferences, alongside real-time market prices. All evaluated strategies successfully fulfill the EV's charging requirements. However, distinct operational patterns are evident: as the risk aversion increases, the EV charging profile becomes smoother and more conservative, thereby avoiding aggressive charging and discharging during volatile price periods and effectively reducing operational risk exposure.

In addition, the proposed method demonstrates excellent computational efficiency. As shown in Table \ref{tab}, the average computational time per decision interval is less than $0.2$ seconds for both risk-neutral and risk-averse strategies, based on a 31-day simulation with five-minute resolution. This efficiency is mainly attributed to the use of closed-form expressions for training post-decision marginal value functions and generating bid curves. Moreover, by directly leveraging price distributions instead of relying on scenario sampling or iterative updates, the method avoids significant computational overhead, making it well suited for real-time implementation.
% \begin{table}
%     \footnotesize
%     \centering
%     \caption{Computational performance of different methods}
%     \label{time}
%     \vspace{-2mm}
%     \begin{tabular}{ccc}
%     \hline \hline
%         Method & \tabincell{c}{offline training \\ (eleven month)} & \tabincell{c}{online optimization \\ (one month)} \\ \hline
%         M1 & $< $  1 mintues & $\approx$ 42.8 mins \\ 
%         M2 & $\backslash$ & $\approx$ 13.6 mins\\ 
%         M3 & $>$ 24 hours & $>$ 24 hours\\ 
%         M4 & $\backslash$ & $>$ 24 hours\\ \hline \hline
%     \end{tabular}
%     \vspace{-1mm}
% \end{table}
% \begin{table}[]
% \begin{tabular}{cccc}
% \hline\hline
%    & Energy cost (\$)  & Degradation cost (\$)  & Total cost (\$)  \\ \hline
% B1 & 8677.83         & 0                    & 8677.83        \\
% B2 & 3442.86         & 530.41               & 3973.27        \\
% B3 & 8082.32         & 115.94               & 8198.26        \\
% S1 & 5763.21         & 184.76               & 5947.97        \\
% S2 & 5517.07         & 207.7                & 5724.77  &   
% \hline\hline
% \end{tabular}
% \end{table}

\section{Conclusion}\label{sec_conclusion}
This paper presents an economic bidding strategy for an EVA participating in real-time electricity markets. By establishing a rigorous analytical connection between electricity price distributions and the marginal opportunity values for the EVA, the proposed method enables interpretable and price-responsive bid construction under both risk-neutral and risk-averse settings. Furthermore, the entire biding process is conducted through closed-form analytical expressions, eliminating the need for scenario-based optimization or numerical solvers. Simulation results using real-world EV charging data from Macao and real-time electricity prices from NYISO validate the proposed strategy, demonstrating superior economic performance, enhanced robustness to price volatility, and high computational efficiency.

\bibliographystyle{ieeetr}
\bibliography{ref}

\begin{thebibliography}{10}

\bibitem{review}
H.~Zhang, X.~Hu, Z.~Hu, and S.~J. Moura, ``Sustainable plug-in electric vehicle integration into power systems,'' {\em Nature Reviews Electrical Engineering}, vol.~1, pp.~35--52, 2024.

\bibitem{ferc}
{Federal Energy Regulatory Commission}, ``Participation of distributed energy resource aggregations in markets operated by rtos and isos (order no.~2222).'' \url{https://www.govinfo.gov/content/pkg/FR-2021-03-30/pdf/2021-06089.pdf}, 2021.
\newblock Federal Register.

\bibitem{ndrc}
{National Development and Reform Commission}, ``Opinions on strengthening the integration and coordinated interaction of new energy vehicles with the power grid.'' \url{https://www.ndrc.gov.cn/xxgk/zcfb/tz/202401/t20240104_1363096.html}, 2024.
\newblock In Chinese.

\bibitem{wuchen}
C.~Lu, J.~Liang, N.~Gu, H.~Wang, and C.~Wu, ``Manipulation-proof virtual bidding mechanism design,'' {\em IEEE Transactions on Energy Markets, Policy and Regulation}, vol.~2, no.~1, pp.~119--131, 2024.

\bibitem{dh1}
Y.~Kabiri-Renani, A.~Arjomandi-Nezhad, M.~Fotuhi-Firuzabad, and M.~Shahidehpour, ``Transactive-based day-ahead electric vehicles charging scheduling,'' {\em IEEE Transactions on Transportation Electrification}, vol.~10, no.~4, pp.~8235--8245, 2024.

\bibitem{dh2}
Z.~Xu, Z.~Hu, Y.~Song, and J.~Wang, ``Risk-averse optimal bidding strategy for demand-side resource aggregators in day-ahead electricity markets under uncertainty,'' {\em IEEE Transactions on Smart Grid}, vol.~8, no.~1, pp.~96--105, 2017.

\bibitem{dh3}
Y.~Chen, Y.~Zheng, S.~Hu, S.~Xie, and Q.~Yang, ``Optimal operation of fast charging station aggregator in uncertain electricity markets considering onsite renewable energy and bounded ev user rationality,'' {\em IEEE Trans lnd. Informat.}, vol.~20, no.~11, pp.~13384--13395, 2024.

\bibitem{rt3}
B.~Luo, Y.~Xu, D.~Xie, Z.~Shi, and H.~Zhang, ``Optimal bidding for aggregated electric vehicles in singapore electricity market: A hierarchical coordination approach,'' {\em IEEE Transactions on Industry Applications}, vol.~61, no.~3, pp.~4913--4923, 2025.

\bibitem{rt2}
Z.~Bao, Z.~Hu, and A.~Mujeeb, ``A novel electric vehicle aggregator bidding method in electricity markets considering the coupling of cross-day charging flexibility,'' {\em IEEE Transactions on Transportation Electrification}, vol.~10, no.~4, pp.~8790--8805, 2024.

\bibitem{rt1}
W.~Wu, J.~Zhu, Y.~Liu, T.~Luo, Z.~Chen, and H.~Dong, ``A coordinated model for multiple electric vehicle aggregators to grid considering imbalanced liability trading,'' {\em IEEE Transactions on Smart Grid}, vol.~15, no.~2, pp.~1876--1890, 2024.

\bibitem{rt4}
C.~Lou, C.~Li, L.~Zhang, W.~Tang, J.~Yang, and J.~Cunningham, ``Two-stage bidding strategy with dispatch potential of electric vehicle aggregators for mitigating three-phase imbalance,'' {\em Journal of Modern Power Systems and Clean Energy}, pp.~1--12, 2025.

\bibitem{xub}
Y.~Baker, N.~Zheng, and B.~Xu, ``Transferable energy storage bidder,'' {\em IEEE Trans. Power Syst.}, vol.~39, no.~2, pp.~4117--4126, 2023.

\bibitem{ecobid}
G.~De~Vivero-Serrano, K.~Bruninx, and E.~Delarue, ``Implications of bid structures on the offering strategies of merchant energy storage systems,'' {\em Applied energy}, vol.~251, p.~113375, 2019.

\bibitem{qiuj}
X.~Lu, J.~Qiu, Y.~Yang, C.~Zhang, J.~Lin, and S.~An, ``Large language model-based bidding behavior agent and market sentiment agent-assisted electricity price prediction,'' {\em IEEE Transactions on Energy Markets, Policy and Regulation}, pp.~1--13, 2024.

\bibitem{kim}
J.~H. Kim and W.~B. Powell, ``An hour-ahead prediction model for heavy-tailed spot prices,'' {\em Energy Economics}, vol.~33, no.~6, pp.~1252--1266, 2011.

\bibitem{drm1}
A.~Ruszczy{\'n}ski, ``Risk-averse dynamic programming for markov decision processes,'' {\em Mathematical programming}, vol.~125, pp.~235--261, 2010.

\bibitem{drm2}
B.~Acciaio and I.~Penner, ``Dynamic risk measures,'' {\em Advanced mathematical methods for finance}, pp.~1--34, 2011.

\bibitem{zhengshuo}
Y.~Li and Z.~Li, ``A novel risk-averse multi-energy management for effective offering strategy of integrated energy production units in a real-time electricity market,'' {\em Applied Energy}, vol.~377, p.~124380, 2025.

\bibitem{zhang1}
H.~Zhang, Z.~Hu, Z.~Xu, and Y.~Song, ``Evaluation of achievable vehicle-to-grid capacity using aggregate pev model,'' {\em IEEE Transactions on Power Systems}, vol.~32, no.~1, pp.~784--794, 2017.

\bibitem{cc1}
J.~Hu, J.~Wu, X.~Ai, and N.~Liu, ``Coordinated energy management of prosumers in a distribution system considering network congestion,'' {\em IEEE Transactions on Smart Grid}, vol.~12, no.~1, pp.~468--478, 2021.

\bibitem{cc2}
M.~González~Vayá and G.~Andersson, ``Self scheduling of plug-in electric vehicle aggregator to provide balancing services for wind power,'' {\em IEEE Trans. Sustain. Energy}, vol.~7, no.~2, pp.~886--899, 2016.

\bibitem{zhang2}
H.~Zhang, Z.~Hu, E.~Munsing, S.~J. Moura, and Y.~Song, ``Data-driven chance-constrained regulation capacity offering for distributed energy resources,'' {\em IEEE Transactions on Smart Grid}, vol.~10, no.~3, pp.~2713--2725, 2019.

\bibitem{powell}
W.~B. Powell, {\em Approximate Dynamic Programming: Solving the curses of dimensionality}, vol.~703.
\newblock John Wiley \& Sons, 2007.

\bibitem{supplydocument}
Z.~Zhu, H.~Zhang, and Y.~Song, ``Supplymentary document.'' \url{https://github.com/Matrix-um/Supplementary-Materials}, 2025.

\bibitem{Lyu}
R.~Lyu, H.~Guo, K.~Zheng, M.~Sun, and Q.~Chen, ``Co-optimizing bidding and power allocation of an ev aggregator providing real-time frequency regulation service,'' {\em IEEE Transactions on Smart Grid}, vol.~14, no.~6, pp.~4594--4606, 2023.

\bibitem{ope}
A.~Philpott, V.~de~Matos, and E.~Finardi, ``On solving multistage stochastic programs with coherent risk measures,'' {\em Operations Research}, vol.~61, no.~4, pp.~957--970, 2013.

\bibitem{math}
V.~Kozm{\'\i}k and D.~P. Morton, ``Evaluating policies in risk-averse multi-stage stochastic programming,'' {\em Mathematical Programming}, vol.~152, pp.~275--300, 2015.

\bibitem{nyiso}
{NYISO}, ``Energy market and operational data.'' \url{https://www.nyiso.com/energy-market-operational-data}, 2025.
\newblock Accessed: 2025-04-01.

\end{thebibliography}

\end{document}